\documentclass[11pt,reqno]{amsart} 

\usepackage[utf8]{inputenc} 
\usepackage{bbm}

\usepackage[margin=1in]{geometry} 

\usepackage{graphicx} 

 \usepackage[parfill]{parskip} 
 
\usepackage{booktabs} 
\usepackage{array} 
\usepackage{paralist} 
\usepackage{verbatim} 
\usepackage{subfig} 
\usepackage{mathrsfs}
\usepackage{amssymb}
\usepackage{amsthm}
\usepackage{amsmath,amsfonts,amssymb,esint}
\usepackage{graphics,color}
\usepackage{enumerate}
\usepackage{mathtools,centernot}
\usepackage{cases}

\usepackage{xfrac}

\newtheorem{theorem}{Theorem}[section]

\newtheorem{corollary}[theorem]{Corollary}

\newtheorem{proposition}[theorem]{Proposition}
\newtheorem{remark}[theorem]{Remark}

\newtheorem*{remark 1}{Remark}

\numberwithin{equation}{section}

\newcommand{\T}{\ensuremath{\mathbb{T}}}
\newcommand*{\R}{\ensuremath{\mathbb{R}}}

\newcommand*{\Z}{\ensuremath{\mathbb{Z}}}

\usepackage{color, graphicx}
\usepackage{mathrsfs, dsfont}

\usepackage{hyperref}

\def\div{\mathop{\rm div}\nolimits}    
\def\curl{\mathop{\rm curl}\nolimits} 

\title{On the Helicity conservation for the incompressible Euler equations}

\author{Luigi De Rosa}
\address{EPFL SB, Station 8, CH-1015 Lausanne, Switzerland}
\email{luigi.derosa@epfl.ch}

\date{\today}
\begin{document}

\begin{abstract}
In this work we investigate the helicity regularity for weak solutions of the incompressible Euler equations. To prove regularity and conservation of the helicity we will threat the velocity $u$ and its $\curl u$ as two independent functions and we mainly show that the helicity is a constant of motion assuming $u \in L^{2r}_t(C^\theta_x)$ and $\curl u \in L^{\kappa}_t(W^{\alpha,1}_x)$ where $r,\kappa $ are conjugate H\"older exponents and $2\theta+\alpha \geq 1$.
Using the same techniques we also show that the helicity has a suitable H\"older regularity even in the range where it is not necessarily constant.
\end{abstract}

\maketitle


\section{Introduction}
In this paper we consider the incompressible Euler equations 
\begin{equation}\label{euler}
\left\{\begin{array}{l}
\partial_t u+ \div (u \otimes u) +\nabla p =0\\ \\
\div u = 0,
\end{array}\right.
\end{equation}
in the spatial periodic setting $\T^3=\R^3\setminus \Z^3$,
where $u: \T^3 \times (0,\infty) \rightarrow \R^3$ is a vector field representing the velocity of the fluid and $p:\T^3 \times (0,\infty) \rightarrow \R$ is the hydrodynamic pressure. 
Letting $\omega := \curl u$, by taking the $\curl$ of the first equation in \eqref{euler} one also gets the evolution equation for the vorticity $\omega$, which is
\begin{equation}\label{vorticity}
\partial_t \omega + \curl \div (u \otimes u) = \partial_t \omega +(u\cdot \nabla) \omega -(\omega\cdot \nabla) u=0\,.
\end{equation}

Thanks to the peculiar structure (and its related cancellation properties) of the non linearity $\div (u\otimes u)$ one can prove that, at least for smooth solutions, we have conservation of quantities like the kinetic energy $E=E(t)$ and the helicity $H=H(t)$. They are defined respectively as
\begin{align*}
E(t)&:= \frac{1}{2}\int_{\T^3} |u|^2(x,t)\, dx \\
H(t)&:=\int_{T^3} u(x,t) \cdot \omega (x,t)\, dx.
\end{align*}

Regarding the kinetic energy it is known that if the solution is sufficiently regular (in space) then it is constant. This was conjectured by the famous physicist Lars Onsager in 1949. He claimed that if $u\in L^\infty_t(C^\theta_x)$, then
\begin{itemize}
\item[$(1)$] for $\theta >\frac{1}{3}$ the kinetic energy is constant;
\item[$(2)$] for $\theta <\frac{1}{3}$ dissipation could appear.
\end{itemize}
Part $(1)$ of the conjecture was completely solved in \cite{CoETi1994}, where the authors proved the energy conservation assuming $u\in L^3_t(B^\theta_{3,\infty})$ for every $ \theta > \frac{1}{3}$ (see also \cite{CCFS2008} for a sharper result). The crucial point of their proof is a careful estimate on the quadratic commutator which arises when one regularize Eq.\eqref{euler} with a standard Friedrichs' mollifier.

The sharpest result in the literature on the helicity conservation has been proved in \cite{CCFS2008} assuming $u \in L^3_t(B^{\sfrac{2}{3}}_{3,c(\mathbb{N})})$. Note that the Sobolev spaces used in this work satisfy $W^{\theta,p}\hookrightarrow B^\theta_{p,c(\mathbb{N})} $, thus one has helicity conservation also for $u \in L^3_t(W^{\frac{2}{3},3}_x)$. Here we propose a different approach which is to threat the velocity and the vorticity as two different functions. We prove the following

\begin{theorem}\label{t:main}
Let $0<\theta,\alpha<1$ and $1 \leq p, q,r,\kappa \leq \infty$ such that $\frac{1}{p}+\frac{1}{q}=\frac{1}{r}+\frac{1}{\kappa}=1$. Suppose that $u$ is a weak solution of \eqref{euler} such that $u \in L^{2r}_t(W^{\theta, 2p}_x)$ and $\omega := \curl u \in L^{\kappa}_t(W^{\alpha,q}_x) $. If $ 2\theta+\alpha\geq 1 $ then $H(t)=H(0)$ for every $t>0$.
\end{theorem}

A similar result to Theorem \ref{t:main} has already been proved in \cite{C2003}. Indeed in \cite{C2003} the author proved the helicity conservation assuming $\omega:=\curl u \in  C^0_t(L^{\sfrac{3}{2}}_x)\cap L^3_t(B^\alpha_{\sfrac{9}{5},\infty})$ for every $\alpha>\frac{1}{3}$. Theorem \ref{t:main} is then a generalization since it treats the velocity and the vorticity separately. Indeed a direct consequence of our theorem is that in order to prove the helicity conservation it suffices to assume $\omega \in L^3_t(W^{\alpha,q}_x)$ for any $\alpha>0$ and any $q>\frac{9}{4+3\alpha}$. We refer to Remark \ref{remark2} for a precise discussion.

Since in our incompressible setting the velocity $u$ is completely determined by its $\curl u$ (thanks to the existence of a potential) then there is a range in which Theorem \ref{t:main} is just a consequence of the conservation proved in \cite{CCFS2008} for $u \in L^3_t(B^{\sfrac{2}{3}}_{3,c(\mathbb{N})})$, and also a range where the hypothesis on $u$ in Theorem \ref{t:main} is redundant.
Thus an interesting case is when the regularity assumption on the $\curl u$ is as weak as possible (see Remark \ref{rem:reg} for a more precise discussion). For this reason one can choose $p=\infty$ and $q=1$ getting the following 

\begin{corollary}\label{c:main}
Let $0<\theta,\alpha<1$ and $1\leq r,\kappa\leq \infty$ such that $\frac{1}{r}+\frac{1}{\kappa}=1$. If  $u\in L^{2r}_t(C^{\theta}_x)$ is a weak solution of \eqref{euler} such that $\omega := \curl u \in L^{\kappa}_t(W^{\alpha,1}_x)$, where $2\theta+\alpha\geq 1$, then the helicity is constant.
\end{corollary}

Note that the hypothesis used in Corollary \ref{c:main} in general do not imply  $u \in L^3_t(B^{\sfrac{2}{3}}_{3,c(\mathbb{N})})$. 

A natural question is to ask whether the helicity as some regularity also in the range in which it is not necessarily constant. To answer this question, instead of the time integrability $L^r_t$ we assume uniformity, namely $L^\infty_t$, showing the following

\begin{theorem}\label{t:second}
Let $0<\theta,\alpha<1$ and $1\leq p, q\leq \infty$ such that $ \frac{1}{p}+\frac{1}{q}=1$. Suppose that $u$ is a weak solution of \eqref{euler} such that $u \in L^{\infty}_t(W^{\theta,2p}_x)$ and $\omega \in L^{\infty}_t(W^{\alpha,q}_x)$. Then there exist a constant $C>0$ such that  
\begin{equation}\label{hold_hel_1}
|H(t)-H(s)|\leq C |t-s|^{\frac{\alpha+\theta}{1-\theta}}\,.
\end{equation}
\end{theorem}

\begin{theorem}\label{t:second2}
Let $\frac{1}{2} < \theta<1$ and suppose that $u$ is a weak solution of \eqref{euler} such that $u \in L^{\infty}_t(W^{\theta,3}_x)$. Then there exist a constant $C>0$ such that  
\begin{equation}\label{hold_hel_2}
|H(t)-H(s)|\leq C |t-s|^{\frac{2\theta-1}{1-\theta}}\,.
\end{equation}
\end{theorem}

We remark that the assumptions $L^\infty_t$ is fundamental in order to get H\"older regularity of $H=H(t)$, but weaker assumptions as $L^r_t$ would also imply suitable Sobolev regularity. However, we are not going to exploit such hypothesis. Moreover the assumption $\theta >\frac{1}{2}$ 

Similar H\"older estimates also hold for the energy $E=E(t)$, see \cite{Is2013}, \cite{CoDe2018}.
The proofs of Theorem \ref{t:second} and Theorem \ref{t:second2} make use of the same techniques introduced in \cite{CoDe2018}, since with this kind of equations one can easily prove H\"older regularity for $E=E(t)$ and $H=H(t)$ by looking at the regularized versions of \eqref{euler} and \eqref{vorticity}.  
Note that the previous theorems still give the helicity conservation if the two H\"older exponents in \eqref{hold_hel_1} and \eqref{hold_hel_2} are bigger than $1$, which means $2\theta+\alpha>1$ and $\theta>\frac{2}{3}$ respectively.  The reader might be confused about the critical hypothesis $2\theta +\alpha=1$ and $\theta=\frac{2}{3}$, which in Theorem \ref{t:second} and Theorem \ref{t:second2} respectively just imply Lipschitz continuity of the helicity instead of conservation, but we remark that the borderline conservation is achieved in Theorem \ref{t:main} and in \cite{CCFS2008} thanks to a limit procedure which is missing in Theorem \ref{t:second}.

Concerning part $(2)$ of the conjecture, in last 10 years an astonishing amount of work has been done in order to produce dissipative H\"older continuous solutions of Euler, see \cite{DLS2012}, \cite{BDLSV2017}, \cite{DLS2015}, \cite{Is2016}. These works are based on a convex integration scheme which has been introduced by C. De Lellis and L. Székelyhidi. For any given (smooth) energy profile $E:[0,T] \rightarrow (0,+\infty)$ and any $\theta<\frac{1}{3}$ these techniques produce solutions $u \in C_{t,x}^\theta$ such that $E(t)=\frac{1}{2}\|u(t)\|_{L^2_x}^2$. Since our Corollary \ref{c:main} shows the conservation of  the helicity if $2\theta +\alpha \geq 1$, then choosing $\theta <\frac{1}{3}$ and the corresponding $\alpha=1-2\theta$, there might exist solutions such that $H=H(t)$ is constant but the energy is not. However we are not able to produce such solutions since in the current works based on convex integration techniques we do not have a strong control on the $\curl u$ in some Sobolev space as the one required here.

Recently these iterative schemes have been adapted also to the hypodissipative Navier-Stokes equations for a sufficiently small power $\gamma$ of the fractional Laplacian $(-\Delta)^\gamma$. Indeed in \cite{CoDLDR2017} and then in \cite{DR2018} the authors proved the existence of infinitely many Leray-Hopf solutions to such equations. The main observation is that (for a small $\gamma$) the dissipative term $(-\Delta)^\gamma$ can be absorbed in the iterative scheme as an error term.

\section{Preliminaries and notations}

\subsection{Helicity and kinetic energy for smooth solutions}\label{helicity_smooth}
Before proving Theorem \ref{t:main} we start considering the helicity for a smooth solution $u$ of \eqref{euler}. By smoothness we can directly compute the first derivative of $H=H(t)$, using equations \eqref{euler} and \eqref{vorticity},  getting

\begin{align*}
\frac{d}{dt}H(t) &= \int_{\T^3} (\partial_t u \cdot \omega + u\cdot \partial_t \omega )\, dx \\
&=-\int_{\T^3} \big( (u\cdot \nabla)u +\nabla p \big) \cdot \omega\, dx -\int_{\T^3} \big( (u\cdot \nabla) \omega-(\omega\cdot \nabla)u\big)\cdot u\, dx \\
&= - \int_{\T^3} \div \big(p \, \omega +u(u \cdot \omega)- \frac{|u|}{2}^2\omega \big)\, dx = 0 \, ,
\end{align*}
where we used the following relations
\begin{align*}
\omega \cdot(u\cdot \nabla) u +u \cdot (u\cdot \nabla)\omega &= \div\big( u(u\cdot \omega) \big) \\
u\cdot (\omega \cdot \nabla) u &= \frac{1}{2} \div \big( |u|^2 \omega \big) \\
\omega \cdot \nabla p &= \div( p \, \omega )\,.
\end{align*}

Thus in the smooth setting, the previous computations easily show that the helicity is constant. 

Similarly, for the kinetic energy, we can multiply the first equation in \eqref{euler} by $u$ getting 
\[
\partial_t \frac{|u|^2}{2}+ u \cdot \div (u \otimes u) +u\cdot \nabla p =\partial_t \frac{|u|^2}{2}+\div \bigg( \frac{|u|^2}{2}u + pu\bigg)=0\,.
\]
Thus, integrating the previous equation over $\T^3$, we can compute
\[
\frac{d}{dt}E(t)=\frac{1}{2}\int_{\T^3}\partial_t|u|^2\,dx= - \int_{\T^3} \div \bigg( \frac{|u|^2}{2}u + pu\bigg)\,dx=0\,.
\]
In order to deal with weak solutions (and so with low regularity), the idea in \cite{CoETi1994} is to mollify the equation \eqref{euler} getting an evolution equation for smooth the quantities $(u_\delta,p_\delta)$, with an "error" forcing therm  which is due to the non-linearity. The crucial observation in \cite{CoETi1994} is that this error has a particular commutator structure and thus satisfies better estimates than $u_\delta$.
Since we also have to deal with the vorticity $\omega$, we will mollify both equations \eqref{euler} and \eqref{vorticity} and we will see that the commutators have exactly the same structure.

\subsection{Spatial H\"older, Sobolev and Besov norms} As already outlined we work in the periodic 3-dimensional spatial domain $\T^3$, thus considering vector fields $u,\omega, f : \T^3 \times (0,+\infty) \rightarrow \R^3$ and a scalar field $p:\T^3 \times (0,+\infty) \rightarrow \R$. We will always denote $\omega := \curl u$.

In what follows $\theta, \alpha \in (0,1)$ and $1\leq p, q < \infty$. We introduce the usual (spatial) H\"older norms as follows. First of all, given any vector field $f$ we will consider its restriction to the $t$-time slice. We thus define the $C_x^0$ norm as
\[ \| f (t)\|_{C^0}:=\sup_{x \in \T^3} \arrowvert f(x,t) \arrowvert\,.
\]
We also define the H\"older seminorm as 
\[
   [f(t)]_{C^\theta} := \sup_{x \neq y, \, x, y \in \T^3}\frac{|f(x,t)-  f(y,t)|}{|x-y|^\theta}\,,
\]
Thus the full H\"older norm is given by
\[
\|f(t)\|_{C^\theta} := \|f(t)\|_0+[f(t)]_\theta\,.
\]
Analogously, we define the  $L^p_x$ and the corresponding $W^{\alpha,p}_x$ norms as
\begin{align*}
\|f(t)\|_{L^p} & := \biggl(\int_{\T^3} |f|^p(x,t)\, dx\biggl)^\frac{1}{p}\,, \\
[f(t)]_{W^{\alpha,p}} & := \biggl(\int_{\T^3} \int_{\T^3} \frac{|f(x,t)-  f(y,t)|^p}{|x-y|^{\alpha p + 3}}\, dx dy\biggl)^\frac{1}{p} \,.
\end{align*}
Then the full Sobolev norm for a fixed time $t$ is given by
\[
\|f(t)\|_{W^{\alpha,p}} := \|f(t)\|_{L^p}+[f(t)]_{W^{\alpha,p}} \, .
\]
For $p=2$ an equivalent norm is given by
$$
[f(t)]_{W^{\alpha,2}}=\| (-\Delta)^{\sfrac{\alpha}{2}}f(t)\|_{L^2}\,.
$$
For $p=\infty$ the Sobolev space can be defined as $W^{\alpha,\infty}\cong C^\alpha$, moreover we set
\[
\|f(t)\|_{L^\infty}:= \sup_{x\in\T^3} |f(x,t)|\,.
\]
We also define the Besov norms as usual
\begin{align*}
[f(t)]_{B^\theta_{p,\infty}}&:= \sup_{y\in \T^3} \frac{\|f(\cdot + y,t)-f(\cdot,t)\|_{L^p}}{|y|^\theta}\, ,\\
\|f(t)\|_{B^\theta_{p,\infty}}&:=\|f(t)\|_{L^p}+[f(t)]_{B^\theta_{p,\infty}}\,.
\end{align*}
For the spaces defined above we have the trivial inclusions $C^\theta  \hookrightarrow B^\theta_{p,\infty}$ and $B^{\theta+\varepsilon}_{p,\infty} \hookrightarrow W^{\theta,p}$ for any $\varepsilon>0$, $1\leq p\leq \infty$, $0<\theta<1$.

In order to avoid confusion, when we have to consider mixed norms in both space and time, we will write explicitly the subscripts $L^p_t$, $W_x^{\alpha,p}$ and $C_x^\theta$. More precisely we will write
\begin{align*}
\|f\|_{L^q_t(W_x^{\alpha,p})} &:= \biggl(    \int_0^{+\infty}  \|f(s)\|_{W^{\alpha,p}} \,ds \biggl)^\frac{1}{q}\,, \\
\|f\|_{L^q_t(C_x^{\theta})} &:= \biggl(    \int_0^{+\infty}  \|f(s)\|_{C^{\theta}} \,ds \biggl)^\frac{1}{q} \,.
\end{align*}

\subsection{Mollification estimates}
Let $ B_1(0)=\{ z \in \mathbb{R}^3 \, : \, |z|<1\}\subseteq \mathbb{R}^3$ be the ball of radius 1 centered in 0  and let  $\rho \in C^ \infty_c (B_1(0))$ such that $\rho \geq 0$ and  $\int_{B_1(0)} \rho(x) dx=1$. For any small parameter $ \delta>0 $ we define the standard convolution kernel by setting $ \rho_\delta:=\delta^{-3} \rho(\frac{x}{\delta})$.
For any function $ f$ we define its mollification (regularization) as 
\[
f_\delta (x):=(f \ast \rho_\delta)(x)=\int_{B_\delta(x)}f(y) \rho_\delta(x-y) \,dy=\int_{B_\delta(0)}f(x-y) \rho_\delta(y)\, dy\,.
\]
It is clear that this definition extends to any $f=f(x,t)$ just taking the space convolution for a fixed time $t$.

In the next proposition we prove some elementary estimates on these regularized functions. We include the proof for the reader convenience. For simplicity we will denote by $\star$ both the scalar and the tensor product between two vectors and for any $1\leq p <\infty$ we set
\[
[ f]_{W^{\theta,p}(\T^3 \llcorner B_\delta)} :=  \biggl(\int_{\T^3} \int_{B_\delta(x)} \frac{|f(x)-  f(y)|^p}{|x-y|^{\alpha p + 3}}\, dx dy\biggl)^\frac{1}{p}
\]

\begin{proposition}\label{moll}
There exists a constant $C$ such that for any $f,g: \mathbb{T}^3 \rightarrow \mathbb{R}^3$ and for any $\theta,\alpha\in (0,1)$ we have:
\begin{align}
\| \nabla f_\delta\|_{L^\infty} &\leq C \delta^{\theta-1} [f]_{C^{\theta}} \, ,\label{est1} \\
 \| f_\delta \star g_\delta -(f \star g)_\delta \|_{L^\infty} &\leq C \delta^{\theta+\alpha} [f]_{C^\theta}[g]_{C^\alpha}\,. \label{est2} 
\end{align}
Moreover for any $1\leq m<\infty$ there exists a positive constant $C=C(m)$ such that for every $1<p, q< \infty$ with $\frac{1}{p}+\frac{1}{q}=1$, we have
\begin{align}
\| \nabla f_\delta\|_{L^m} &\leq C \delta^{\alpha-1} [f]_{W^{\alpha,m}(\T^3 \llcorner B_\delta)} \, ,\label{est4} \\
\| \nabla \curl f_\delta\|_{L^m} &\leq C \delta^{\alpha-2} [f]_{W^{\alpha,m}(\T^3 \llcorner B_\delta)} \, ,\label{est7} \\
 \| f_\delta \star g_\delta -(f \star g)_\delta \|_{L^m} &\leq C \delta^{\theta+\alpha} [f]_{W^{\theta,mp}(\T^3 \llcorner B_\delta)}[g]_{W^{\alpha,mq}(\T^3 \llcorner B_\delta)}\,. \label{est5}\\
 \| f_\delta \star g_\delta -(f \star g)_\delta \|_{L^m} &\leq C \delta^{\theta+\alpha} [f]_{C^\theta}[g]_{W^{\alpha,m}(\T^3 \llcorner B_\delta)}\,. \label{est6}
\end{align}
\end{proposition}

Note that the previous proposition is stated for time-independent functions $f,g$, thus applying it for a fixed $t$-time slice we get the same estimates for any time dependent vector field.
\begin{proof}

Since for any $\delta>0$ we have $ \| \nabla \rho_\delta \|_{C^0} \leq C \delta^{-4}$, for some constant $C$ which depends only on $\nabla \rho$ in $B_1(0)$, we can estimate
 \begin{align}\label{estproof1}
 | \nabla f_\delta|(x) &=  \bigg|\int_{B_\delta (0)} f(x-y) \otimes \nabla \rho_\delta(y) \, dy \bigg|   =  \bigg|\int_{B_\delta (0)} \big( f(x-y)-f(x) \big) \otimes \nabla \rho_\delta(y) \, dy \bigg|  \nonumber \\
 &\leq C \delta^{-4}  \int_{B_\delta (0)} \big| f(x-y)-f(x) \big|  \, dy   \leq C \delta^{-4+\theta} \int_{B_\delta (0)} \frac{\big| f(x-y)-f(x) \big|}{|y|^{\theta}}  \, dy\,.
 \end{align}
Taking the $L^\infty$ norm on both sides, estimating  $\int_{B_\delta (0)} \frac{\big| f(x-y)-f(x) \big|}{|y|^{\theta}}\, dy \leq  C \delta^3 [f]_{C^\theta} $, we get \eqref{est1}. Notice also that taking the power $p$ of \eqref{estproof1} and using Jensen's inequality we achieve
\[
| \nabla f_\delta|^p(x)\leq C \delta^{(\theta -1)p} \int_{B_\delta (0)} \frac{\big| f(x-y)-f(x) \big|^p}{|y|^{\theta p}}  \, \frac{dy}{\delta^3} \leq C \delta^{(\theta -1)p} \int_{B_\delta (0)} \frac{\big| f(x-y)-f(x) \big|^p}{|y|^{\theta p+3}} \, dy \, ,
\]
from which, integrating over $\T^3$ and taking the $p$-root, one gets \eqref{est4}. In the same way one can prove \eqref{est7} putting the operator $\nabla \curl$ on the kernel $\rho_\delta$ and this gives an extra $\delta^{-1}$.

 Now for every $ x \in \T^3$, since  $\int_{B_\delta(0)} \rho_\delta(y) \,dy=1 $, we have
\begin{align}
(f \star g)_\delta(x) -f_\delta(x) \star g_\delta(x)&=\int_{B_\delta(0)} \big(f(x-y)-f(x) \big)\star \big(g(x-y)-g(x) \big) \rho_\delta(y) \, dy \nonumber \\
& - \int_{B_\delta(0)} \big(f(x-y)-f(x) \big)\rho_\delta(y) \, dy\star \int_{B_\delta(0)} \big(g(x-y)-g(x) \big)\rho_\delta(y) \, dy ,\nonumber
\end{align}
and again, since $ |y| \leq \delta$ we get 
\begin{align}\label{estproof2}
\big|(f \star g)_\delta -f_\delta \star g_\delta \big|(x)& \leq C\delta^{\theta+\alpha} \int_{B_\delta(0)} \frac{\big|f(x-y)-f(x) \big|}{|y|^{\theta}} \frac{\big|g(x-y)-g(x) \big|}{|y|^{\alpha}}  \rho_\delta(y) \, dy \nonumber \\
&+C \delta^{\theta+\alpha}  \int_{B_\delta(0)}\frac{\big|f(x-y)-f(x) \big|}{|y|^\theta}\rho_\delta(y) \, dy \int_{B_\delta(0)}\frac{\big|g(x-y)-g(x) \big|}{|y|^\alpha}\rho_\delta(y) \, dy \nonumber \\
& \leq C \delta^{\theta+\alpha}[f]_{C^\theta}[g]_{C^\alpha} \, ,
\end{align} 
 which  proves \eqref{est2}. We now conclude with the proof of \eqref{est5} (estimate \eqref{est6} is then easier and is left to the reader). We observe that for any $x\in \T^3$, using Jensen and H\"older inequalities, we have
 \begin{align}\label{primo}
 &\bigg| \int_{B_\delta(0)} \big(f(x-y)-f(x) \big)\star \big(g(x-y)-g(x) \big) \rho_\delta(y) \, dy\bigg|^m \nonumber \\
 &\leq C \int_{B_\delta(0)} \big|f(x-y)-f(x) \big|^m \big|g(x-y)-g(x) \big|^m\rho_\delta(y) \, dy\nonumber\\
 & \leq C \bigg( \int_{B_\delta(0)} |f(x-y)-f(x)|^{pm}\rho_\delta(y)\,dy \bigg)^{\frac{1}{p}} \bigg(  \int_{B_\delta(0)} |g(x-y)-g(x)|^{qm}\rho_\delta(y)\,dy \bigg)^{\frac{1}{q}} \nonumber\\
 &\leq C \delta^{m(\theta+\alpha)}\bigg( \int_{B_\delta(0)} \frac{|f(x-y)-f(x)|^{pm}}{|y|^{3+mp\theta}}\,dy \bigg)^{\frac{1}{p}} \bigg(  \int_{B_\delta(0)}  \frac{|g(x-y)-g(x)|^{pm}}{|y|^{3+mp\alpha}}\,dy \bigg)^{\frac{1}{q}}\,.
 \end{align}
 Similarly we estimate
  \begin{align}\label{secondo}
 &\int_{\T^3} \bigg| \int_{B_\delta(0)} \big(f(x-y)-f(x) \big) \rho_\delta(y) \, dy \star\int_{B_\delta(0)} \big(g(x-y)-g(x) \big) \rho_\delta(y) \, dy \bigg|^m \,dx\nonumber \\
 &\leq C \int_{\T^3} \bigg(  \int_{B_\delta(0)} \big|f(x-y)-f(x) \big|^m \rho_\delta(y) \, dy\int_{B_\delta(0)} \big|g(x-y)-g(x) \big|^m \rho_\delta(y) \, dy \bigg) \,dx\nonumber\\
 & \leq C \bigg( \int_{\T^3}\int_{B_\delta(0)} |f(x-y)-f(x)|^{pm}\rho_\delta(y)\,dy \bigg)^{\frac{1}{p}} \bigg(  \int_{\T^3}\int_{B_\delta(0)} |g(x-y)-g(x)|^{qm}\rho_\delta(y)\,dy \bigg)^{\frac{1}{q}} \nonumber\\
 &\leq C \delta^{m(\theta+\alpha)}\bigg( \int_{\T^3}\int_{B_\delta(0)} \frac{|f(x-y)-f(x)|^{pm}}{|y|^{3+mp\theta}}\,dy \bigg)^{\frac{1}{p}} \bigg(  \int_{\T^3} \int_{B_\delta(0)}  \frac{|g(x-y)-g(x)|^{pm}}{|y|^{3+mp\alpha}}\,dy \bigg)^{\frac{1}{q}}\,.
 \end{align}
 Finally putting together \eqref{primo} and \eqref{secondo} and using once again H\"older inequality with $p, q$ we get \eqref{est6}.
 \end{proof}
 
 \section{Proofs of the main theorems}\label{main_section}
Before proving our main results we start with two remarks about the hypothesis needed in Theorem \ref{t:main}.
\begin{remark}\label{rem:reg}
Theorem \ref{t:main} is stated for any couple of exponents $1 \leq p,q\leq \infty$ such that $\frac{1}{p}+\frac{1}{q}=1$. The reader may wonder if the assumption on $u$ could be redundant, since it could be a consequence of the one on the $\curl u$. Indeed in our incompressible setting we have $u=\curl \big((-\Delta)^{-1} \curl u \big)$.
In particular, if $\curl u \in W^{\alpha,q}_x$, $1<q<\infty$, by Calderón–Zygmund we get $u \in W_x^{1+\alpha,q}$ and by Sobolev embeddings we have that $W^{1+\alpha,q} \hookrightarrow W^{\theta,\frac{2q}{1-q}}$ if
\begin{equation}\label{embedding}
q>\frac{9}{5+2(\alpha-\theta)}\,.
\end{equation}
In the case $q=1$ we have  $u \in W_x^{1+\alpha-\varepsilon,1}$ for any $\varepsilon>0$, but this is obviously not enough to guarantee any H\"older regularity on $u$.
\end{remark}

\begin{remark}\label{remark2}
For any  $\alpha>0$ we assume $\omega\in W^{\alpha,q}_x$ and we choose $\theta=\frac{1-\alpha}{2}$, so that the helicity is preserved. Then by \eqref{embedding} we have that $u \in W^{\frac{1-\alpha}{2},\frac{2q}{1-q}}_x$ if 
$$
q>\frac{9}{4+3\alpha}\,.
$$
Thus we have that the assumption $\alpha>\frac{1}{3}$ in \cite{C2003} is not necessary if one assume more on the integrability exponent $q$.
\end{remark}

We highlight that if $u\in W^{\frac{1}{2},2}_x$ then the helicity is the action of the functional $\curl u\in  W^{-\frac{1}{2},2}_x$ on the velocity $u$ and it can be rapresented as
\begin{equation}\label{hel_dual}
H(t)=\int_{\T^3}(-\Delta)^{\frac{1}{4}} u\cdot  (-\Delta)^{-\frac{1}{4}} \curl u \,dx\,.
\end{equation}
Note that by Cauchy-Schwarz and Calderón–Zygmund we have
$$
|H(t)|\leq \|(-\Delta)^{\frac{1}{4}} u (t)\|_{L^2{(\T^3)}} \|(-\Delta)^{-\frac{1}{4}} \curl u(t)\|_{L^2{(\T^3)}} \leq C \|u(t)\|_{W^{\frac{1}{2},2}(\T^3)}^2\,.
$$
\subsection{Proof of Theorem \ref{t:main}}We first mollify equations \eqref{euler} and \eqref{vorticity} getting
\begin{align}
\partial_t u_\delta+ \div (u_\delta \otimes u_\delta) +\nabla p_\delta &=\div R_\delta\, , \label{euler_moll} \\
\partial_t \omega_\delta +(u_\delta\cdot \nabla) \omega_\delta -(\omega_\delta \cdot \nabla) u_\delta &= \curl \div R_\delta\, , \label{vorticity_moll}
\end{align}
where $R_\delta:= u_\delta \otimes u_\delta- (u \otimes u)_\delta$. Now we consider the helicity $H_\delta$ related to the smooth vector fields $u_\delta, \omega_\delta$, namely the function 
\begin{equation}\label{helmoll}
H_\delta(t) := \int_{\T^3} u_\delta(x,t) \cdot  \omega_\delta(x,t) \,dx\,.
\end{equation}
By the regularity of $u$ and $\omega$ it is clear that for almost every $t\geq 0$, $ H_\delta(t)\rightarrow H(t)$ as $\delta \rightarrow 0$. We can now compute the time derivative of $H_\delta$ (in order to be precise at this point one should also mollify $u$ and $\omega$ in time, say with some $\rho_\varepsilon=\rho_\varepsilon(t)$, in order to rigorously  reach inequality \eqref{ineq} for $u_{\delta,\varepsilon}$ and $\omega_{\delta,\varepsilon}$ and conclude that \eqref{ineq} holds letting $\varepsilon \rightarrow 0$). Using \eqref{euler_moll} and \eqref{vorticity_moll} as in Section \ref{helicity_smooth}, we get
\begin{align}\label{der:hel}
\frac{d}{dt}H_\delta(t) &= - \int_{\T^3} \div \big(p_\delta \, \omega_\delta +u_\delta(u_\delta \cdot \omega_\delta)- \frac{|u_\delta|}{2}^2\omega_\delta \big)\, dx \nonumber \\ 
& + \int_{\T^3} \omega_\delta \cdot \div R_\delta \, dx + \int_{\T^3} u_\delta \cdot \curl \div R_\delta \, dx \nonumber \\
& = -2\int_{\T^3}\nabla \omega_\delta : R_\delta \, dx\, ,
\end{align}
where in the last equality we used the integration by parts formula. Thus we have that
\begin{align}\label{ineq}
\big| H_\delta(t)- H_\delta(0)\big|\leq 2 \int_0^t \int_{\T^3} |\nabla \omega_\delta|(x,s)  |R_\delta|(x,s) \, dx ds \leq 2\int_0^t \| \nabla \omega_\delta(s) \|_{L^q} \| R_\delta(s) \|_{L^p}\,ds \,,
\end{align}
and by Proposition \ref{moll} we conclude that 
\[
\big| H_\delta(t)- H_\delta(0)\big|\leq  C \delta^{2\theta +\alpha-1} \int_0^t [\omega(s)]_{W^{\alpha,q}( \T^3 \llcorner B_\delta)} [u(s)]^2_{W^{\theta,2p}( \T^3 \llcorner B_\delta)}\, ds\,.
\]
Note that in the previous estimate we have used two conjugate exponents $1<p,q<\infty$, the case where one of them is equal to $1$ (or equivalently $\infty$) is analogous. Finally, using H\"older inequality with exponents $r,\kappa$ we achieve 
\[
\big| H_\delta(t)- H_\delta(0)\big|\leq  C \delta^{2\theta +\alpha-1} [\omega]_{L^\kappa_t(W_x^{\alpha,q}( \T^3 \llcorner B_\delta))} [u]^2_{L^{2r}_t(W_x^{\theta,2p}( \T^3 \llcorner B_\delta))}\, ,
\]
thus the claim follows letting $\delta \rightarrow 0$.

\subsection{Proof of Theorem \ref{t:second}}
Now we will see how the $L^\infty_t$ assumption leads to some H\"older regularity of the helicity even without the assumption $2\theta +\alpha \geq 1$. The following technique comes from \cite{CoDe2018}, where the authors proved H\"older regularity for the kinetic energy assuming $u\in L^\infty_t(C^\theta _x)$.

We define $H_\delta(t)$ as in \eqref{helmoll}. For any couple of times $s,t$ we estimate
\begin{equation}\label{step1}
\big| H(t)-H(s)\big| \leq \big| H(t)-H_\delta(t)\big|+\big| H_\delta(t)-H_\delta(s)\big|+\big| H_\delta(s)-H(s)\big|\,.
\end{equation}

By the $L^\infty_t$ assumption, both the first and the third term can be estimated using \eqref{est5} with $m=1$ and $p=q=2$ as follows
\[
\big| H(t)-H_\delta(t)\big|+\big| H_\delta(s)-H(s)\big| \leq C \delta^{\theta+\alpha} \|u\|_{L^\infty_t(W^{\theta,p}_x)} \|\omega\|_{L^\infty_t(W^{\alpha,q}_x)}\,,
\]
where, in order to apply \eqref{est5}, we have also used the the property $H(t)=\int_{\T^3} u\cdot \omega=\int_{\T^3} (u\cdot \omega)_\delta$.
We are left with the second summand in the right hand side of \eqref{step1}. We have that
\[
\big| H_\delta(t)-H_\delta(s)\big|\leq |t-s| \bigg\| \frac{d}{dt}H_\delta\bigg\|_{L^\infty_t}\, ,
\]
and by \eqref{der:hel} together with Proposition \ref{moll} we get 
\[
\big| H_\delta(t)-H_\delta(s)\big|\leq C |t-s|  \delta^{2\theta +\alpha-1}\|u\|^2_{L^\infty_t(W^{\theta,2p}_x)} \|\omega\|_{L^\infty_t(W^{\alpha,q}_x)}\,.
\]
Combining the previous estimates with \eqref{step1} we have achieved
\[
\big| H(t)-H(s)\big| \leq C  (\delta^{\theta +\alpha}+ |t-s|\delta^{2\theta +\alpha-1})\, ,
\]
for some constant $C>0$ which depends on both $u, \omega$. Finally choosing $\delta = |t-s|^{\frac{1}{1-\theta}}$ we can conclude
\[
\big| H(t)-H(s)\big| \leq C|t-s|^{\frac{\alpha+\theta}{1-\theta}}\,.
\]
\subsection{Proof of Theorem \ref{t:second2}}
The proof runs in the same way as the one for Theorem \ref{t:second}. By equation \eqref{der:hel} and using \eqref{est7} and \eqref{est5} we have
$$
\big| H_\delta(t)-H_\delta(s)\big|\leq |t-s| \bigg\| \frac{d}{dt}H_\delta\bigg\|_{L^\infty_t}\leq 2 |t-s| \| \nabla \curl u_\delta\|_{L^\infty_t(L^3_x)} \|R_\delta\|_{L^{\infty}_t(L^{\sfrac{3}{2}}_x)}\leq C |t-s| \delta^{3\theta-2}[ u]^3_{L^\infty_t(W^{\theta,3}_x)}\,. 
$$
Since for every $\delta>0$
$$
H(t)=\int_{\T^3} (-\Delta)^{\sfrac{1}{4}} u\cdot  (-\Delta)^{-\sfrac{1}{4}}  \curl u\,dx=\int_{\T^3} \big( (-\Delta)^{\sfrac{1}{4}} u \cdot (-\Delta)^{-\sfrac{1}{4}} \curl u)_\delta\,dx\,,
$$
applying \eqref{est5} with $m=1$ we deduce that for every $t\geq 0$
$$
\big| H(t)-H_\delta(t)\big|\leq C \delta^{2\theta-1} [(-\Delta)^{\sfrac{1}{4}}  u]_{L^\infty_t(W_x^{\theta-\sfrac{1}{2},2})} [(-\Delta)^{-\sfrac{1}{4}} \curl u]_{L^\infty_t(W_x^{\theta-\sfrac{1}{2},2})}\leq C\delta^{2\theta-1}[u]^2_{L^\infty_t(W^{\theta,2}_x)}\,,
$$
where in the last inequality we also used Calderón–Zygmund estimates.
Thus we achieved
\begin{align*}
\big| H(t)-H(s)\big|& \leq\big| H(t)-H_\delta(t)\big| +\big| H_\delta(t)-H_\delta(s)\big|+\big| H_\delta(s)-H(s)\big| \\
&\leq  C  \Big(\delta^{2\theta -1}[u]^2_{L^\infty_t(W^{\theta,2}_x)}+ |t-s|\delta^{3\theta -2}[u]^3_{L^\infty_t(W^{\theta,2}_x)}\Big)\, ,
\end{align*}
from which we conclude by choosing $\delta=|t-s|^{\frac{1}{1-\theta}}$.

\end{document}